\documentclass[fleqn,a4paper,12pt]{article}
\usepackage{amssymb,amsfonts,amsmath,amsthm}

\numberwithin{equation}{section}

\begin{document}


\begin{center}
\Large{
\textbf{Intersections of Magnus subgroups and embedding theorems for cyclically presented groups}
}

\bigskip

\normalsize{
Martin Edjvet
}

School of Mathematical Sciences\\
University of Nottingham\\
University Park\\
Nottingham NG7 2RD, UK

\bigskip

James Howie

Department of Mathematics\\
Heriot-Watt University\\
Edinburgh EH14 4AS, UK
\end{center}


\section{Introduction}

Let $G=\langle Y | r \rangle$ be a one-relator group where $r$ is a cyclically reduced word in the
free group on $Y$.  A subset $Y_1 \subseteq Y$ is called a \textit{Magnus subset} if $Y_1$ omits a
generator which appears in the relator $r$.  A subgroup $M_1$ of $G$ is called a \textit{Magnus subgroup}
if $M_1 = \langle Y_1 \rangle$ for some Magnus subset $Y_1$ of $Y$ and so by the Magnus Freiheitssatz
\cite{M} (or \cite[page 109]{LS})
$M_1$ is free of rank $|Y_1|$.  There have been recent studies in 
\cite{C,Ho}  on intersections of such subgroups.  In fact it is shown in \cite{C} that the intersection of two Magnus
subgroups $M_i = \langle Y_i \rangle$ ($1 \leq i \leq 2$) is either the free group $F(Y_1 \cap Y_2)$
on $Y_1 \cap Y_2$ or is the free product of $F(Y_1 \cap Y_2)$ together with an infinite cyclic group.
When the latter holds we say that $M_1$ and $M_2$ have \textit{exceptional intersection}.

Suppose now that $F_n$ denotes the free group of rank $n$ generated by the elements of the set
$X=\{ x_0,\ldots ,x_{n-1} \}$ and let $\theta \colon F_n \to F_n$ be the
automorphism of $F_n$ for which $x_i \theta = x_{i+1}$ (where subscripts are taken
mod $n$).  Let $w$ be a cyclically reduced element of $F_n$ and define the group
$G_n(w)=F_n /N$ where $N$ is the normal closure in $F_n$ of the set
$\{ w,w \theta , \ldots , w \theta ^{n-1} \}$.  Then a group $G$ is said to have a
\textit{cyclic presentation} or is \textit{cyclically presented} if $G \cong G_n (w)$ for
some $n$ and for some $w$.

In recent papers \cite{EHT,E,HR} there has been an interest in finding which cyclically
presented groups define the trivial group.  On the other hand the papers 
\cite{EH,N,P} 
primarily give examples of infinite cyclically presented groups.

In this paper we will use results from \cite{C,Ho} to give further examples of infinite cyclically
presented groups.  In particular we will prove two embedding theorems.

Before stating our first theorem we will need some further notation and definitions.  Assume
that $k\geq 1$ and that $n \geq 3k$.  Let $F_{3k}$ denote the subgroup of $F_n$ of
rank $3k$ generated by the subset $\{ x_0,\ldots ,x_{3k-1}\}$ of $X$; and let $F_{k+1}$
denote the subgroup of $F_n$ of rank $k+1$ generated by the subset $\{ x_0,\ldots ,x_k \}$
of $X$.  Let $w=w(x_0,\ldots ,x_k)$ be a cyclically reduced element of $F_{k+1}$ that
involves both $x_0$ and $x_k$ (and perhaps other generators).  Define the groups $G$ and
$H_1$ by putting $G=G_n(w)$ and $H_1=F_{3k}/K$ where $K$ is the normal closure in
$F_{3k}$ of the set $\{ w,w\theta ,\ldots ,w \theta^{2k-1} \}$.  Thus we have
\[
G=\langle x_0,\ldots ,x_{n-1} \mid w(x_0,\ldots ,x_k),\ldots ,w(x_{n-1},\ldots ,x_{k-1})\rangle
\]
and
\[
H_1=\langle x_0,\ldots ,x_{3k-1} \mid w(x_0,\ldots ,x_k), \ldots , w(x_{2k-1},\ldots,
x_{3k-1})\rangle .
\]
Given this we can now state the following.

\textbf{Theorem 1.1} \quad \textit{If
$\langle x_0,\ldots ,x_{k-1} \rangle \cap \langle x_1,\ldots ,x_k \rangle \cap \ldots \cap
\langle x_k, \ldots , x_{2k-1} \rangle = \{1\}$ in the group $H_1$ then
$\langle x_i,\dots,x_{i+k} \mid w \theta^i \rangle$ embeds in $G=G_n(w)$ for
$0 \leq i \leq n-1$. In particular, $G$ is an infinite group.}

We prove this result in Section 2.  It is worth noting here the following corollary to Theorem 1.1.

 \textbf{Corollary 1.2} \quad
\textit{If the Magnus subgroups
$\langle x_0 ,\ldots , x_{k-1} \rangle$ and $\langle x_1,\ldots ,x_k \rangle$ do not have
exceptional intersection in the one-relator group
$\langle x_0,\ldots ,x_k \mid w \rangle$ where $w=w(x_0,\ldots ,x_k)$ involves both
$x_0$ and $x_k$ then $G_n(w)$ is infinite for $n \geq 4k$.}

The proof of Corollary 1.2 is as follows.
If we regard the group $H_1$ of Theorem 1.1 as a stem product of
one-relator groups then the condition in the statement of Theorem 1.1 is a condition about
the intersections of Magnus subgroups.  Since by assumption these intersections are
non-exceptional we have
$\langle x_0,x_1,\ldots ,x_{k-1} \rangle \cap \langle x_1,\ldots ,x_k \rangle \cap \ldots \cap
\langle x_k,\ldots ,x_{2k-1} \rangle$
$=\langle x_1,\ldots ,x_{k-1} \rangle \cap \langle x_2,\ldots,x_{k+1} \rangle \cap \ldots \cap
\langle x_k,\ldots ,x_{2k-1} \rangle$ $\ldots$
$=\langle x_{k-1} \rangle \cap \langle x_k,\ldots ,x_{2k-1} \rangle = \{1\}$ and the result
follows.

It is interesting to reflect on the bound $n \geq 4k$ in Theorem 1.1 and Corollary 1.2.  Using
different methods it may be possible to improve upon this.  However as things stand
we must have $n\geq 3k+1$
since for example $H_{3k}=G_{3k} (x_0^{-1} x_k x_0 x_k^{-2} )$ is trivial for $k \geq 1$ \cite{Hi}.
Observe that there is a `gap' of $k$ between the two subscripts involved in $w$ in $H_{3k}$.
The influence of this gap between subscripts seems significant and this is reflected in our next result
which will allow for an improvement on $4k$ (and indeed $3k+1$)
provided there are no exceptional intersections and the maximum gap that occurs is not too big.

First we need the following definition.
Let $0<t \leq k$.  Then a cyclically reduced word $w$ in the alphabet $\{ x_0,\dots,x_k \}$ will be
called $t$-\textit{pure} if:
\begin{enumerate}
\item[1.]
$w$ involves $x_0$ and $x_k$;
\item[2.]
the Magnus subgroups $M_0 = \langle x_1,\dots,x_k \rangle$ and
$M_k = \langle x_0, \dots, x_{k-1} \rangle$ of $$H=\langle x_0,\dots,x_k \mid w \rangle$$ do not
have exceptional intersection;
\item[3.]
for each $i=1,\dots,k-t$, at least one of the letters $x_i,\dots,x_{i+t-1}$ is involved in $w$, so that
the subsets $X_{i,i+t-1}=\{ x_0,\dots,x_{i-1}, x_{i+t},\dots,x_k \}$ of $X$ freely generate
Magnus subgroups $M_{i,i+t-1}$ of $H$; and
\item[4.]
for each $i=1,\dots,k-t$, the Magnus subgroup $M_{i,i+t-1}$ does not have exceptional intersection
with either $M_0$ or $M_k$.
\end{enumerate}

Given this we can now state

\textbf{Theorem 1.3} \quad \textit{Let $w$ be a $t$-pure, cyclically reduced word in
$\{ x_0,\dots,x_k \}$.  Then for $i=0,\dots,n-1$ the group
$\langle x_i,\dots,x_{i+k} \mid w \theta^i \rangle$
embeds in $G_n (w)$ for $n \geq 2k+2t$. In particular $G_n (w)$ is infinite.}

The following corollary is immediate.

\textbf{Corollary 1.4} \quad \textit{If the Magnus subgroups $\langle x_0,\ldots,x_{k-1} \rangle$
and $\langle x_1,\ldots,x_k \rangle$ do not have exceptional intersection in the one-relator group
$\langle x_0,\ldots,x_k \mid \omega \rangle$ where $\omega = \omega (x_0,\ldots,x_k)$
involves $x_i$ for $0 \leq i \leq k$ then $G_n (\omega)$ is infinite for $n \geq 2(k+1)$.}

We will prove Theorem 1.3 in Section 3 and we end the paper with some further remarks and some
examples in Section 4.

\section{Proof of Theorem 1.1}

We adapt an argument due originally to Higman \cite{Hi}.
In addition to the groups $G$ and $H_1$ defined in the introduction we will also need the
group $H_2$ defined by
\[
H_2 = \langle x_{2k},\ldots , x_{n-1},x_0,\ldots ,x_{k-1} \mid
w(x_{2k},\ldots , x_{3k}), \ldots , w(x_{n-1},\ldots ,x_{k-1})\rangle .
\]

\textbf{Proof of Theorem 1.1} \quad Let $A_1=\langle x_0,\ldots ,x_{k-1} \rangle$ and
$B_1=\langle x_{2k},\ldots ,x_{3k-1} \rangle$ be subgroups of $H_1$; and let
$A_2=\langle x_{2k},\ldots ,x_{3k-1} \rangle$ and
$B_2=\langle x_0,\ldots ,x_{k-1} \rangle$ be subgroups of $H_2$.  Now $H_1$ can be
expressed as a stem product of one-relator groups
\begin{align*}
\langle x_0,\ldots ,x_k \mid w(x_0,&\ldots ,x_k)\rangle \underset{C_1}{\ast}
\langle x_1,\ldots ,x_{k+1} \mid w(x_1,\ldots ,x_{k+1})\rangle \underset{C_2}{\ast}
\ldots\\
&\ldots \underset{C_{2k-1}}{\ast}
\langle x_{2k-1},\ldots ,x_{3k-1} \mid w(x_{2k-1},\ldots ,x_{3k-1})\rangle
\end{align*}
where $C_j=\langle x_j,\ldots ,x_{j+k-1} \rangle$ for $1\leq j\leq 2k-1$; and a similar
statement holds for $H_2$.  It follows that $A_1$ and $B_1$ are free subgroups of $H_1$
of rank $k$ and that $A_2$ and $B_2$ are free subgroups of $H_2$ of rank $k$.

It is now sufficient to prove the claim that $\langle A_1,B_1 \rangle = A_1 \ast B_1$ in $H_1$.
To see that this is sufficient, the fact that $n\geq 4k$ will allow us to deduce in a similar way
that $\langle A_2,B_2 \rangle = A_2 \ast B_2$ in $H_2$.  Thus
$\langle A_1,B_1 \rangle \leq H_1$ and $\langle A_2,B_2 \rangle \leq H_2$ are free subgroups
each of rank $2k$ and $G$ is the amalgamated free product
($H_1 \ast H_2$; $\langle A_1,B_1 \rangle = \langle A_2,B_2 \rangle$) where we are
identifying $x_i$ in $\langle A_1,B_1 \rangle$ with $x_i$ in $\langle A_2,B_2 \rangle$ for
each $i \in \{ 0,\ldots ,k-1,2k,\ldots ,3k-1 \}$.
But each $\langle x_i,\dots, x_{i+k} \mid w \theta^i \rangle$ embeds in either $H_1$ or $H_2$
for $0 \leq i \leq n-1$ and so embeds in $G$.

To establish the claim observe that the stem product $H_1$ acts (without inversion) on a tree
$T$ with fundamental region a tree with $2k$ vertices $u_i$ ($0 \leq i \leq 2k-1$) and $2k-1$ edges
$e_j$ ($0 \leq j \leq 2k-2$) where $e_j$ joins the vertex $u_j$ to $u_{j+1}$ for
$j=0,\dots,2k-2$, where the stabilizer $V_j= \text{Stab}_{H_1} (u_j)$ of the vertex $u_j$ is the
one-relator group $\langle x_j,\ldots ,x_{j+k} \mid w(x_j,\ldots ,x_{j+k}) \rangle$ for
$0\leq j \leq 2k-1$ and the stabilizer $E_i= \text{Stab}_{H_1}(e_i)$ of the edge $e_i$ is
the subgroup $\langle x_{i+1},\ldots ,x_{i+k} \rangle$ of $V_i \cap V_{i+1}$.  Since
$\{ x_{i+1},\ldots ,x_{i+k} \}$ omits $x_i$ and $x_{i+k+1}$ it follows from the
Magnus Freiheitssatz \cite{M} that $E_i$ is a free group of rank $k$.  (For the basic theory of
groups acting on trees see \cite{B} or \cite{Se}.)

Let $a_1 \in A_1 \setminus \{1\}$ and suppose that $a_1 \in E_{k-1}$.  Then $a_1$ fixes
the vertices $u_0$ and $u_k$ and so belongs to the stabilizer of the geodesic in $T_1$ from
$u_0$ to $u_k$.  It follows that $a_1 \in E_j$ for $0\leq j \leq k-1$ and so
$a_1 \in \langle x_0,\ldots ,x_{k-1} \rangle \cap \ldots \cap \langle x_k,\ldots ,x_{2k-1}
\rangle = \{1\}$ by the assumption in the statement of the theorem.  This contradiction
shows that $A_1 \cap E_{k-1} = \{1\}$.  Since $T$ is a tree we can by deleting the edge
$e_{k-1}$ partition the vertex set $V(T)$ of $T$ into the disjoint union $Y_1 \overset{\cdot}{\cup}
Y_2$ in which $u_{k-1} \in Y_1$ and $u_k \in Y_2$ (so in particular vertex $v \in Y_1$ if and
only if there is a path in $T$ from $v$ to $u_{k-1}$ which omits $e_{k-1}$).  A
distance-preserving argument now shows that $y_2 a_1 \in Y_1$ ($\forall y_2 \in Y_2$)
($\forall a_1 \in A_1 \setminus \{ 1 \}$).  We can now apply symmetry to deduce that
$B_1 \cap E_{k-1} = \{1\}$ and that $y_1 b_1 \in Y_2$ ($\forall y_1 \in Y_1$)
($\forall b_1 \in B_1 \setminus \{1\}$).  Finally $\langle A_1,B_1 \rangle = A_1 \ast B_1$
now follows immediately by a Ping-Pong argument: if $w=a_1b_1\cdots a_mb_m$
is a cyclically reduced word in $A_1\ast B_1$ with $m\ge 1$, then 
$u_0a_1=u_0\in Y_1$, so $u_0a_1b_1\in Y_2$,
$u_0a_1b_1a_2\in Y_1$ and so on. Hence $u_0w\in Y_2$, so $u_0w\ne u_0$, and so $w\ne 1$ in $H_1$.
  $\square$

\section{Proof of Theorem 1.3}

We first prove the following lemma on graphs of groups.  For the basic definitions and theory, we refer the reader to,
for example, \cite{B}. We use the convention that all graphs are oriented; the 
initial and terminal vertices of an edge $e$ are denoted $\iota(e),\tau(e)$
respectively.

\textbf{Lemma 3.1} \quad \textit{Let $\Gamma$ be a tree, and $(\mathcal{G},\Gamma)$ and
$(\mathcal{M},\Gamma)$ be graphs of groups.  Denote the vertex groups by $G_v$ and $M_v$,
the edge groups $G_e$ and $M_e$, and the fundamental groups by $G$, $M$ respectively.}
\textit{Suppose that, for each $v$, we have an injective homomorphism $\phi_v \colon M_v \to G_v$
such that}
\[
\phi_{\iota (e)}^{-1} (G_e)=M_e=\phi_{\tau (e)}^{-1} (G_e) \quad
\mathrm{and} \quad \phi_{\iota(e)}|_{M_e}=\phi_{\tau(e)}|_{M_e} \quad (\mathrm{for~each~edge}~e~\mathrm{of}~\Gamma).
\]
\textit{Then these induce an injective homomorphism $\phi \colon M \to G$.}

\textit{Proof}.  For each edge $e$ of $\Gamma$, the maps $M_e \to M_{\iota (e)} \to
G_{\iota (e)} \to G$ and $M_e \to M_{\tau (e)} \to G_{\tau (e)} \to G$ defined by
$\phi_{\iota (e)}$ and $\phi _{\tau (e)}$ respectively coincide, by hypothesis.  Hence the $\phi_v$
do indeed induce a homomorphism $\phi \colon M \to G$.  It remains to show that $\phi$ is injective.

Now the groups $M$ and $G$ act on trees $T_M$, $T_G$, respectively, with quotient $\Gamma$
in each case.  Moreover, the trees $T_M$, $T_G$ each contain a copy of $\Gamma$ as a fundamental
domain for the action.  More specifically, there are injections $i_M \colon \Gamma \to T_M$ and
$i_G \colon \Gamma \to T_G$ such that $M_x$ is the stabiliser in $M$ of
$i_M (x) \in T_M$ and $G_x$ is the stabiliser in $G$ of $i_G (x) \in T_G$,
for each vertex or edge $x$ of $\Gamma$.  The group $M$ also acts
on the tree $T_G$ via the homomorphism $\phi \colon M \to G$.  Moreover, the isomorphism
$i_G \circ i_M^{-1} \colon i_M(\Gamma) \to i_G (\Gamma)$ between the fundamental domains
extends uniquely to an $M$-equivariant graph-map $\Phi \colon T_M \to T_G$ defined on vertices and edges by
\[
\Phi \Big( m \big( i_M (v) \big) \Big) = \phi (m) \big( i_G (v) \big), \quad
\Phi \Big( m \big( i_M (e) \big) \Big) = \phi (m) \big( i_G (e) \big)
\quad \mathrm{respectively}.
\]
By hypothesis, if $f$ is an edge of $T_M$ with one vertex in $i_M (\Gamma)$ and
$\Phi (f)=i_G(e) \in i_G(\Gamma)$, then the fact that $f=m \big( i_M (e) \big)$ for some
$m \in M$ can be used to show that $f=i_M (e)$.  To see this suppose, for example,
$\iota (f) = i_M (v)$ for $v=\iota (e) \in \Gamma$.  Then
$m \in \text{Stab}_M \big( i_M (v) \big) = M_v$ and
$\phi (m) \in \text{Stab}_G \big( i_G (e) \big) = G_e$, so
$m \in \phi_v^{-1} (G_e)=M_e=\text{Stab}_M \big( i_M (e) \big)$ and so
$f=m \big( i_M (e) \big) = i_M (e)$.  A similar argument holds if
$\tau (f) \in i_M (\Gamma)$.  This shows that $\Phi$ is locally injective at vertices of
$i_M (\Gamma)$.  Since $i_M (\Gamma)$ is a fundamental domain for the action of $M$ on $T_M$
and $\Phi$ is $M$-equivariant, it follows that $\Phi$ is locally injective at all vertices of $T_m$,
that is, an immersion $T_M \to T_G$.  Since $T_G$ is a tree, $\Phi \colon T_M \to T_G$ is injective.

Now suppose that $m \in \textrm{Ker} (\phi)$ and $v$ is a vertex of $\Gamma$.  Then
$\Phi \Big( m \big( i_M (v) \big) \Big) = \phi (m) \big( i_G (v) \big)=i_G (v) = \Phi \big( i_M (v) \big)$.
Since $\Phi$ is injective, $m \big( i_M (v) \big) = i_M (v)$, whence $m \in M_v$.  But by hypothesis
$\phi |_{M_v}=\phi_v \colon M_v \to G_v$ is injective, so $m=1$.

Hence $\phi$ is injective, as claimed.  $\Box$

\textbf{Proposition 3.2} \quad \textit{Let $w$ be a $t$-pure, cyclically reduced word in the alphabet
$\{ x_0,\dots, x_k \}$.  Then the set $A=\{ x_0,\ldots,x_{k-1},x_{k+t}, \ldots, x_{2k+t-1}\}$
freely generates a free subgroup of}
\[
\hat{G}=\langle x_0, \ldots, x_{2k+t-1} \mid w, w \theta , \ldots, w \theta^{k+t-1} \rangle .
\]
\textit{Proof}.  Let $\Gamma$ be the tree with vertices $v_0, \ldots, v_{k+t-1}$ and edges $e_i$
joining $v_{i-1}$ to $v_i$ for $i=1,\ldots, k+t-1$.  Then $\hat{G}=\pi_1 (\mathcal{G},\Gamma)$,
where the vertex groups are $G_i = G_{v_i} = \langle x_i,\ldots,x_{i+k} \mid w \theta^i \rangle$
and the edge groups are the free groups $G_{e_i}=\langle x_{i+1}, \ldots, x_{i+k} \rangle$.
Note that, by hypothesis, these are Magnus subgroups embedded in the adjacent vertex groups.

Define $M_i = M_{v_i}$ to be the free group on $A \cap \{ x_i, \ldots, x_{i+k} \}$, and
$M_{e_i}$ to be the free group on $A \cap \{ x_{i+1}, \ldots, x_{i+k} \}$ for each $i$.  Then this
defines a graph of groups $(\mathcal{M},\Gamma)$ whose fundamental group
$M= \pi_1 ( \mathcal{M},\Gamma)$ is free on $A$.

Define maps $\phi_i \colon M_i \to G_i$ by inclusion of generating sets, and note that these maps
satisfy the hypotheses of Lemma 3.1 by the $t$-pure condition.  For example, that $\phi_i$ is
injective follows from the fact that $w$ cannot omit $t$ consecutive generators from the list
$\{ x_0,\ldots, x_k \}$ and involves $x_0,x_k$, while the fact that $\phi_i^{-1} (G_{e_i})=
M_{e_i}=\phi_{i+1}^{-1} (G_{e_{i+1}})$ follows from the non-exceptionality of the corresponding Magnus
intersections. Finally, the fact that $\phi_i$ and $\phi_{i+1}$
agree on $M_{e_i}$ is immediate from the fact that they agree on
generators.

It follows that the induced map $M \to \hat{G}$ is injective, so $A$ freely generates a
subgroup of $\hat{G}$, as claimed. $\Box$

The proof of Theorem 1.3 is similar to that of Theorem 1.1.  We have
$$G=G_n(w) = \langle x_0,\ldots,x_{n-1} \mid w(x_0,\ldots,x_k), \ldots, w(x_{n-1},\ldots,
x_{k-1}) \rangle.$$  
This time we put
$H_1 = \langle x_0,\ldots,x_{2k+t-1} \mid w(x_0,\ldots,x_k), \ldots, w(x_{k+t-1}, \ldots,
x_{2k+t-1}) \rangle$ and
$H_2 = \langle x_{k+t}, \ldots, x_{n-1}, x_0, \ldots, x_{k-1} \mid w(x_{k+t}, \ldots,
x_{2k+t}), \ldots, w(x_{n-1}, x_0, \ldots, x_{k-1}) \rangle$.  It follows from Proposition 3.2
and the fact that $n \geq 2k + 2t$ that the set $$\{ x_0, \ldots, x_{k-1}, x_{k+t}, \ldots,
x_{2k+t-1}\}$$ freely generates a free subgroup in both $H_1$ and $H_2$ so that as in the proof
of Theorem 1.1 the group $G$ can be expressed as an amalgamated free product of $H_1$ and
$H_2$ and the result follows.

\section{Concluding remarks}

\subsection{}
Observe that in Corollary 1.2 and 1.4 non-exceptional intersection is required for only two specific
Magnus subgroups.  For example if we consider
$$\langle x_0,x_1,x_2 \mid (x_0 x_1)^{-1} x_2 (x_0 x_1)x_2^{-2} \rangle$$ then
$x_0^{-1} x_2 x_0 = x_1 x_2^2 x_1^{-1}$ and so
$\langle x_0,x_2 \rangle \cap \langle x_1,x_2 \rangle=\langle x_2,x_0^{-1} x_2 x_0 \rangle=
\langle x_2,x_1 x_2^2 x_1^{-1} \rangle$ has exceptional intersection.  However
$\langle x_0,x_1 \rangle \cap \langle x_1,x_2 \rangle = \langle x_1 \rangle$ is
non-exceptional and so $G_n \big( (x_0 x_1)^{-1} x_2 (x_0 x_1) x_2^{-1} \big)$ is
infinite for $n\geq 6$ by Corollary 1.4.

\subsection{}
The following consequence of applying Corollary 4.2 and Lemma 5.1 of 
\cite{Ho} to our situation
provides a combinatorial method for checking for the presence of exceptional intersections.

\textbf{Proposition 4.1} \quad \textit{The intersection
$\langle x_0,\ldots ,x_{k-1} \rangle \cap \langle x_1,\ldots ,x_k \rangle$ in
$$\langle x_0,\ldots ,x_k \mid w(x_0,\ldots ,x_k) \rangle$$ is exceptional only if
$w(x_0,\ldots ,x_k)$ is of one of the following two forms:}
\begin{enumerate}
\item[(i)]
\textit{$w_1^{\alpha _1} w_2^{\beta _1} w_1^{\alpha _2} w_2^{\beta _2} \ldots
w_1^{\alpha _l} w_2^{\beta _l}$ where $w_1 \in \langle x_0,\ldots ,x_{k-1} \rangle$
and $w_2 \in \langle x_1,\ldots ,x_k \rangle$; or}
\item[(ii)]
\textit{$w_3^{\alpha _1} (v_1 v_2)^{\beta _1} w_3^{\alpha _2} (v_1 v_2)^{\beta _2}
\ldots w_3^{\alpha _l} (v_1 v_2)^{\beta _l}$ where $w_3 \in \langle x_1,\ldots ,x_{k-1}
\rangle$, $v_1 \in \langle x_0,\ldots ,x_{k-1} \rangle$ and
$v_2 \in \langle x_1,\ldots ,x_k \rangle$;}
\end{enumerate}
\textit{where $\alpha _j,\beta _j \in \mathbb{Z}$ for $1\leq j\leq l$.}

\subsection{}
When there are exceptional intersections in Corollary 1.2 it could happen that
Theorem 1.1 may still apply.  For example let
$w=(x_0 x_2)^{-1} x_1 (x_0 x_2)x_1^{-2}$.  Then $w=1$ implies
$x_0^{-1} x_1 x_0 = x_2 x_1^2 x_2^{-1}$ and so
$\langle x_0,x_1 \rangle \cap \langle x_1,x_2 \rangle = \langle x_1,x_2 x_1^2 x_2^{-1}
\rangle$; and again $\langle x_1,x_2 \rangle \cap \langle x_2,x_3 \rangle =
\langle x_2,x_1^{-1} x_2 x_1 \rangle$.  But
$\langle x_1,x_2 x_1^2 x_2^{-1} \rangle \cap \langle x_2,x_1^{-1} x_2x_1 \rangle=\{1\}$
and it follows that $\langle x_0,x_1 \rangle \cap \langle x_1,x_2 \rangle \cap
\langle x_2,x_3 \rangle = \{1\}$ and $G_n(w)$ is infinite for $n\geq 8$.

\subsection{}
If $\langle x_0,\ldots ,x_{k-1} \rangle \cap \langle x_1,\ldots ,x_k \rangle=
\langle x_0,\ldots ,x_{k-1} \rangle$ as for the groups
$G_n (x_0[x_1^\alpha ,x_2^\beta ])$ studied in \cite{EH} then
$\langle x_0,\ldots ,x_{k-1} \rangle \leq \langle x_j,\ldots ,x_{j+k-1} \rangle$
($1\leq j \leq k$) and Theorem 1.1 does not apply.  Perhaps the graph-immersion approach
of Gersten \cite{G} and Stallings \cite{St} may yield information on when exactly
$\langle x_0,\ldots ,x_{k-1} \rangle \cap \ldots \cap \langle x_k,\ldots ,x_{2k-1} \rangle
\neq \{1\}$.


\begin{thebibliography}{99}
\bibitem{B} 
G Baumslag, \textit{Topics in Combinatorial Group Theory}, Lectures in Mathematics ETH
Z\"urich (Birkh\"auser) (1993).
\bibitem{C} 
D J Collins, Intersections of Magnus subgroups of one-relator groups, \textit{Groups: topological,
combinatorial and arithmetic aspects}, LMS Lecture Notes 311 (CUP) (2004) 255--296.
\bibitem{EHT}
M Edjvet, P Hammond and N Thomas, Cyclic presentations of the trivial group,
\textit{Experiment.\ Math.\ }\textbf{10} (2001) 303--306.
\bibitem{E} 
M Edjvet, On cyclic presentations, \textit{J.~Group Theory} \textbf{6} (2003) 261--270.
\bibitem{EH}
M Edjvet and P Hammond, On a class of cyclically presented groups,
\textit{International J.~Algebra and Computation} \textbf{14} (2004) 213--240.
\bibitem{G} 
S M Gersten, Intersections of finitely generated subgroups of free groups and resolutions of
graphs, \textit{Inventiones Mathematicae} \textbf{71} (1983) 567--591.
\bibitem{HR}
G Havas and E F Robertson, Irreducible cyclic presentations of the trivial group,
\textit{Experiment.\ Math.\ }\textbf{12} (2003) 487--490.
\bibitem{Hi}
G Higman, A finitely generated infinite simple group, \textit{J.~London Math.~Soc.~}\textbf{26}
(1951) 61--64.
\bibitem{Ho}  
J Howie, Magnus intersections in one-relator products, \textit{Michigan Math.~J.~}\textbf{53}
(2005) 597--623.
\bibitem{LS}
R C Lyndon and P E Schupp, \textit{Combinatorial Group Theory}
(Springer) (1977).
\bibitem{M}
W Magnus,
\"{U}ber diskontinuierliche Gruppen mit einer definierenden Relation
(Der Freiheitssatz), \textit{J. reine angew. Math.} \textbf{163} (1930) 141-165.
\bibitem{N} 
B H Neumann, Some group presentations, \textit{Canadian J.~Math.~}\textbf{30} (1978)
838--850.
\bibitem{P} 
S J Pride, Groups with presentations in which each relator involves exactly two generators,
\textit{J.~London Math.\ Soc.\ }\textbf{36} (1987) 245--256.
\bibitem{Se} 
J-P Serre, \textit{Trees} (translated from the French by J Stillwell) (Springer) (1980).
\bibitem{St} 
J R Stallings, Topology of finite graphs, \textit{Inventiones Mathematicae} \textbf{71}
(1983) 551--565.
\end{thebibliography}
\end{document}